\documentclass[12pt,twoside]{amsart}

\usepackage{amsfonts}
\usepackage{xr}
\usepackage{graphicx}
\usepackage{amsmath}
\usepackage{epsfig}
\usepackage{color}
\def\R{\mathbb{R}}
\def\N{\mathbb{N}}
\def\C{\mathbb{C}}
\def\epsilon{\varepsilon}

\def\tilde{\widetilde}

\def\trait (#1) (#2) (#3){\vrule width #1pt height #2pt depth #3pt}
\def\fin{\hfill\trait (0.1) (5) (0) \trait (5) (0.1) (0) \kern-5pt 
\trait (5) (5) (-4.9) \trait (0.1) (5) (0)}

\newcommand{\be}{\begin{equation}}
\newcommand{\ee}{\end{equation}}
\newcommand{\baa}{\begin{array}}
\newcommand{\eaa}{\end{array}}
\newcommand{\ba}{\begin{eqnarray}}
\newcommand{\ea}{\end{eqnarray}}

\newcommand{\Bk}{\color{black}}

\newtheorem{theo}{\bf Theorem}[section]
\newtheorem{lem}[theo]{\bf Lemma}
\newtheorem{pro}[theo]{\bf Proposition}
\newtheorem{cor}[theo]{\bf Corollary}

\newtheorem{rem}[theo]{\bf Remark}

\title[Fractional Poincar\'e inequalities]{\bf{Fractional Poincar\'e
    inequalities for general measures}} \author{Cl\'ement Mouhot}
\author{Emmanuel Russ} \author{Yannick Sire}

\begin{document}

\begin{abstract} We prove a fractional version of Poincar\'e
  inequalities in the context of $\R^n$ endowed with a fairly general
  measure. Namely we prove a control of an $L^2$ norm by a non local
  quantity, which plays the role of the gradient in the standard
  Poincar\'e inequality. The assumption on the measure is the fact
  that it satisfies the classical Poincar\'e inequality, so that our
  result is an improvement of the latter inequality. Moreover we also
  quantify the tightness at infinity provided by the control on the
  fractional derivative in terms of a weight growing at infinity. The
  proof goes through the introduction of the generator of the
  Ornstein-Uhlenbeck semigroup and some careful estimates of its
  powers. To our knowledge this is the first proof of fractional
  Poincar\'e inequality for measures more general than L\'evy
  measures.
\end{abstract} 

\maketitle 

\tableofcontents

\section{Introduction}

The aim of this paper is to prove a Poincar\'e inequality on $\R^n$,
endowed with a measure $M(x) \, dx$, involving nonlocal quantities in
the right-hand side in the spirit of Gagliardo semi-norms for Sobolev
spaces $W^{s,p}(\R^n)$ (see e.g. \cite{adams}).

Throughout this paper, we denote by $M$ a positive weight in
$L^1(\R^n)$. We assume that $M$ is a $C^2$ function and that this
measure $M$ satisfies the usual Poincar\'e inequality: there exists a
constant $\lambda(M) > 0$ such that $\forall \, f \in H^1(\R^n, M), $
\begin{equation}
 \label{eq:poincare} 
 \int_{\R^n} |\nabla f(y)|^2 \, M(y) \, dy \ge 
 \lambda(M) \, \int_{\R^n} \left| f(y) - \int_{\R^n} f(x) \, M(x) \,
   dx \right|^2 \, M(y) \, dy.
\end{equation}
If the measure $M$ can be written $M=e^{-V}$, this inequality is known
to hold (see \cite{BBCG}, or also \cite{villaniams}, Appendix A.19,
Theorem 1.2, see also \cite{deuschstroock}, Proof of Theorem 6.2.21
for related criteria) whenever there exist $a \in (0,1)$, $c >0$ and
$R>0$ such that 
\begin{equation} \label{assumptionV}
\forall \, |x| \ge R, \quad 
a \, |\nabla V(x)|^2 - \Delta V \ge c.
\end{equation}
In particular, the inequality (\ref{eq:poincare}) holds, for instance,
when $M=(2\pi)^{-n/2}\exp\left(-\left\vert x\right\vert^2/2\right)$ is
the Gaussian measure, but also when $M(x)=e^{-\left\vert
    x\right\vert}$, and more generally when $M(x)=e^{-\left\vert
    x\right\vert^{\alpha}}$ with $\alpha\ge 1$. Note that, when $V$ is
convex and
$$
\mbox{Hess}(V)\geq \mbox{cst} \, \mbox{Id}
$$
on the set where $\left\vert V\right\vert<+\infty$, the measure
$M(x)dx$ satisfies the log-Sobolev inequality, which in turn implies
(\ref{eq:poincare}) (see \cite{ledouxlivre}).

In the sequel, by $L^2(\R^n,M)$, we mean the space of measurable
functions on $\R^n$ which are square integrable with respect to the
measure $M(x) \, dx$, by $L^2_0(\R^n,M)$ the subspace of functions of
$L^2(\R^n,M)$ such that $\int_{\R^n} f(x) \, M(x) \, dx=0$, and by
$H^1(\R^n,M)$, the Sobolev space of functions in $L^2(\R^n,M)$, the
weak derivative of which belongs to $L^2(\R^n,M)$.

As it shall be proved to be useful later on, remark that, under a
slightly stronger assumption than (\ref{assumptionV}), the Poincar\'e
inequality \eqref{eq:poincare} admits the following self-improvement:
\begin{pro} \label{poincimprovedM} 
 Assume that $M$ there exists $\varepsilon>0$ such that
\begin{equation} \label{alpha}
\frac{(1-\varepsilon)\left\vert \nabla V\right\vert^2}{2} 
-\Delta V \xrightarrow[x \to \infty]{}
  +\infty
, \qquad M = e^{-V}.
\end{equation}
Then there exists $\lambda'(M)>0$ such that, for all function $f\in
L^2_0(\R^n,M) \cap H^1(\R^n,M)$:
 \begin{equation} \label{pim} 
   \iint_{\R^n} \left\vert \nabla f(x)\right\vert^2 \, M(x) \, dx 
   \ge \lambda'(M) \, \int_{\R^n} \left\vert f(x)\right\vert
   ^2\left(1+\left\vert \nabla \ln M(x)\right\vert^2\right) \, M(x) \, dx.
\end{equation}
\end{pro}

The proof of Proposition \ref{poincimprovedM} is classical and will be
given in Appendix A for sake of completeness. \Bk \bigskip

\noindent We want to generalize the inequality \eqref{eq:poincare} in
the strenghtened form of Proposition \ref{poincimprovedM}, replacing,
in the right-hand side, the $H^1$ semi-norm by a non-local expresssion
in the flavour of the Gagliardo semi-norms.

We establish the following theorem:

\begin{theo} \label{mainth} Assume that $M=e^{-V}$ is a $C^2$ positive
  $L^1$ function which satisfies \eqref{alpha}. Let $\alpha\in
  (0,2)$. Then there exist $\lambda_\alpha(M)>0$ and $\delta(M)$
  (constructive from our proof and the usual Poincar\'e constant
  $\lambda'(M)$) such that, for any function $f$ belonging to a dense
  subspace of $L^2_0(\R^n,M)$ , we have
  \begin{eqnarray}\label{poincfrac}
    \iint _{\R^n\times \R^n}
    \frac{\left\vert f(x)-f(y)\right\vert^2}{\left\vert
        x-y\right\vert^{n+\alpha}} \, M(x) \, e^{-\delta(M)\Bk \, |x-y|} \Bk \, dx
    \, dy \ge \\ \nonumber
    \lambda_\alpha(M) \, \int_{\R^n} \left\vert
      f(x)\right\vert^2 \big (1+\left\vert \nabla \ln
      M(x)\right\vert^{\alpha}\big ) \, M(x) \, dx.
        \end{eqnarray} 
\end{theo}

\begin{rem} 
  Inequality~\eqref{poincfrac} could as usual be extended to any
  function $f$ with zero average such that both sides of the
  inequality make sense.  In particular it is satisfied for any
  function $f$ with zero average belonging to the domain of the
  operator $L=-\Delta - \nabla V \cdot \nabla $ that we shall
  introduce later on.  As we shall see, this domain is dense in
  $L^2_0(\R^n,M)$.
\end{rem}
\medskip

Observe 
 that the left-hand side of (\ref{poincfrac}) involves a fractional
 moment of order $\alpha$ related to the homogeneity of the semi-norm
 appearing in the right-hand side. One could expect in the left-hand
 side of (\ref{poincfrac}) the Gagliardo semi-norm for the fractional
 Sobolev space $H^{\alpha/2}(\R^n,M)$, namely
 $$
 \iint _{\R^n\times \R^n} \frac{\left\vert
     f(x)-f(y)\right\vert^2}{\left\vert x-y\right\vert^{n+\alpha}}
 M(x) \, M(y) \, dx \, dy.
 $$
 Notice that, instead of this semi-norm, we obtain a ``non-symmetric''
 expression. However, our norm is more natural: one should think of
 the measure over $y$ as the L\'evy measure, and the measure over $x$
 as the ambient measure. We emphasize on the fact that our measure is
 rather general and in particular, as a corollary of Theorem
 \ref{mainth}, we obtain an automatic improvement of the Poincar\'e
 inequality \eqref{eq:poincare} by
 
 \begin{eqnarray*}
   \iint _{\R^n\times \R^n}
    \frac{\left\vert f(x)-f(y)\right\vert^2}{\left\vert
        x-y\right\vert^{n+\alpha}} \, M(x)  \, dx
    \, dy \ge 
        \lambda_\alpha(M) \, \int_{\R^n} \left\vert
      f(x)\right\vert^2   \, M(x) \, dx.
        \end{eqnarray*} 

The  question  of  obtaining  Poincar\'e-type  inequalities  (or  more
generally entropy  inequalities) for  L\'evy operators was  studied in
the probability  community in  the last decades.  For instance  it was
proved by Wu \cite{wu} and Chafa\"{\i} \cite{chafai} that
\[
\mbox{Ent}_\mu ^\Phi (f) \le \int \Phi''(f) \, \nabla f \cdot \sigma
\cdot \nabla f \, d\mu + \int \int D_\Phi(f(x),f(x+z)) \, d\nu_\mu(z)
\, d\mu(x)
\]
(see also the use of this inequality in \cite{GI}) with 
\[
\mbox{Ent}_\mu ^\Phi (f) = \int \Phi(f) \, d\mu - \Phi\left( \int f \, d\mu
\right)
\]
and $D_\Phi$ is the so-called Bregman distance associated to $\Phi$:
\[
D_\Phi(a,b) = \Phi(a) - \Phi(b) - \Phi'(b)\, (a-b),
\]
where $\Phi$ is some well-suited functional with convexity properties,
$\sigma$ the matrix of diffusion of the process, $\mu$ a rather
general measure, and $\nu_\mu$ the (singular) L\'evy measure
associated to $\mu$.  Choosing $\Phi(x)=x^2$ and $\sigma=0$ yields a
Poincar\'e inequality for this choice of measure $(\mu,\nu_\mu)$. The
improvement of our approach is that we do not impose any link between
our measure $M$ on $x$ and the singular measure $|z|^{-n-\alpha}$ on
$z=x-y$. This is to our knowledge the first result that gets rid of
this strong constraint.

\begin{rem}
  Note that the exponentially decaying factor $e^{-\delta(M)\Bk \,
    |x-y|}$ in \eqref{poincfrac} also {\em improves} the inequality as
  compared to what is expected from Poincar\'e inequality for L\'evy
  measures. This decay on the diagonal could most probably be further
  improved, as shall be studied in futur works. Other extensions in
  progress are to allow more general singularities than the Martin
  Riesz kernel $\frac{1}{|x-y|^{n+\alpha}}$ (see the book
  \cite{landkof}) and to develop an $L^p$ theory of the previous
  inequalities.
\end{rem}

Our proof heavily relies on fractional powers of a (suitable
generalization of the) Ornstein-Uhlenbeck operator, which is defined
by
$$
Lf=-M^{-1}\, \mbox{div} (M\, \nabla f)=-\Delta f-\nabla \ln
M\cdot\nabla f,
$$
for all $f\in {\mathcal D}(L):=\left\{g\in H^1(\R^n,M);\
  \mbox{div}(M\nabla g)\in L^2(\R^n)\right\}$.  One therefore has, for
all $f\in {\mathcal D}(L)$ and $g\in H^1(\R^n,M)$,
$$
\int_{\R^n} Lf(x)g(x)M(x) \, dx=\int_{\R^n} \nabla f(x)\cdot \nabla
g(x) \, M(x) \, dx.
$$
It is obvious that $L$ is symmetric and non-negative on $L^2(\R^n,M)$,
which allows to define the usual power $L^{\beta}$ for any $\beta\in
(0,1)$ by means of spectral theory.  Note that $L^{\alpha/2}$ is {\em
  not} the symmetric operator associated to the Dirichlet form $\iint
_{\R^n\times \R^n} \frac{\left\vert f(x)-f(y)\right\vert^2}{\left\vert
    x-y\right\vert^{n+\alpha}} \, M(x) \, dx \, dy$.

We now describe the strategy of our proofs. The proof of Theorem
\ref{mainth} goes in three steps. We first establish $L^2$
off-diagonal estimates of Gaffney type on the resolvent of $L$ on
$L^2(\R^n,M)$. These estimates are needed in our context since we do
not have Gaussian pointwise estimates on the kernel of the operator
$L$.
 
Then, we bound the quantity
$$
\int_{\R^n} \left\vert f(x)\right\vert^2\left(1+\left\vert
    \nabla\ln M(x)\right\vert^{\alpha}\right) \, M(x) \, dx
$$ 
in terms of $\left\Vert L^{\alpha/4}f\right\Vert_{L^2(\R^n,M)}
^2$. This will be obtained by an abstract argument of functional
calculus based on rewriting in a suitable way the conclusion of
Proposition \ref{poincimprovedM}. Finally, using the $L^2$
off-diagonal estimates for the kernel of $L$, we establish that
$$
\left\Vert L^{\alpha/4}f\right\Vert_{L^2(\R^n,M)}^2 \leq C\iint
_{\R^n\times \R^n} \frac{\left\vert f(x)-f(y)\right\vert^2}{\left\vert
    x-y\right\vert^{n+\alpha}} \, M(x) \, dx \, dy,
$$
which would conclude the proof. 
\medskip

As can be seen from the rough sketch previously described, we borrow
methods from harmonic analysis. This seems not so common in the field
of Poincar\'e and Log-Sobolev inequalities (to the knowledge of the
authors), where standard techniques rely on global functional
inequalities, see for instance the powerful so-called
$\Gamma_2$-calculus of Bakry and Emery (\cite{be}). We hope this paper
will stimulate further exchanges between these two fields.

\section{Off-diagonal $L^2$ estimates for the resolvent of $L$}

We recall that for every $f \in {\mathcal D}(L)$, we define
\begin{equation}
 \label{eq:defL}
 Lf=-M^{-1} \, \mbox{div} (M \, \nabla f)=-\Delta f-\nabla \ln M\cdot\nabla f
\end{equation}
 From the fact that $L$ is self-adjoint and nonnegative on $L^2(\R^n,M)$ we have
\[ 
\| (L-\mu)^{-1} \|_{L^2(\R^n,M)} \le
\frac{1}{\mbox{dist}(\mu,\Sigma(L))}
\]
where $\Sigma(L)$ denotes the spectrum of $L$, and $\mu \not \in
\Sigma(L)$.  Then we deduce that $(\mbox{I}+t \, L)^{-1}$ is bounded with
norm less than $1$ for all $t >0$. Since $tL(\mbox{I}+tL)^{-1}=\mbox{I}-(\mbox{I}+tL)^{-1}$, the same is true for $tL(\mbox{I}+tL)^{-1}=\mbox{I}-(\mbox{I}+tL)^{-1}$ with a norm less than $2$. Moreover, $\nabla(\mbox{I}+tL)^{-1}\in H^{1}(\R^n,M)$.\par
\noindent Actually, when $f\in L^2(\R^n,M)$ is supported in a closed set $E\subset \R^n$ and $F\subset \R^n$ is a closed subset disjoint from $E$, a much more precise estimate on the $L^2$ norm of $(\mbox{I}+t \, L)^{-1}f$ and $tL(\mbox{I}+t \, L)^{-1}f$on $F$ can be given. Here are these $L^2$ off-diagonal estimates for the
resolvent of $L$:\Bk 
\begin{lem} \label{off} There exists $C_1=C_1(M)>0$ (constructive from our
  proof) with the following property: for all closed disjoint subsets
  $E,F\subset \R^n$ with $\mbox{dist}(E,F)=:d>0$, all function $f\in
  L^2(\R^n,M)$ supported in $E$ and all $t>0$,
  $$
  \left\Vert (\mbox{I}+t \, L)^{-1}f\right\Vert_{L^2(F,M)}+\left\Vert t \,
    L(\mbox{I}+t \, L)^{-1}f\right\Vert_{L^2(F,M)}\leq 8 \, e^{-C_1 \,
    \frac{d}{\sqrt t}} \left\Vert f\right\Vert_{L^2(E,M)}.
  $$
\end{lem}

Note that, in different contexts, this kind of estimate, originating in
\cite{gaff}, turns out to be a powerful tool, especially when no pointwise
upper estimate on the kernel of the semigroup generated by $L$ is available
(see for instance \cite{auscherams, kato, amr}).
Since we found no reference for these off-diagonal estimates for the
resolvent of $L$, we give here a proof. \par

\medskip

\noindent {\bf Proof of Lemma \ref{off}.} We argue as in \cite{kato}, Lemma
1.1. From the fact that $L$ is self-adjoint on $L^2(\R^n,M)$ we have
\[ 
\| (L-\mu)^{-1} \|_{L^2(\R^n,M)} \le
\frac{1}{\mbox{dist}(\mu,\Sigma(L))}
\]
where $\Sigma(L)$ denotes the spectrum of $L$, and $\mu \not \in
\Sigma(L)$.  Then we deduce that $(\mbox{I}+t \, L)^{-1}$ is bounded with
norm less than $1$ for all $t >0$, and it is clearly enough to argue when
$0<t<d$.

Define $u_t=(\mbox{I}+t \, L)^{-1}f$, so that, for all function $v\in
H^{1}(\R^n,M)$,
\begin{eqnarray} \label{test} \int_{\R^n} u_t(x) \, v(x) \, M(x) \,
 dx+ t \, \int_{\R^n} \nabla u_t(x)\cdot \nabla v(x) \, M(x) \,
 dx=\\  \nonumber
  \int_{\R^n} f(x) \, v(x) \, M(x) \, dx.
\end{eqnarray}
Fix now a nonnegative function $\eta\in {\mathcal D}(\R^n)$ vanishing
on $E$. Since $f$ is supported in $E$, applying (\ref{test}) with
$v=\eta^2 \, u_t$ (remember that $u_t\in H^1(\R^n,M)$) \Bk yields
\[
\int_{\R^n} \eta^2(x)\left\vert u_t(x)\right\vert^2
\, M(x) \, dx + t \, \int_{\R^n} \nabla u_t(x)\cdot \nabla
(\eta^2u_t) \, M(x) \, dx=0,
\]
which implies
\begin{multline*}
 \int_{\R^n} \eta^2(x)\left\vert u_t(x)\right\vert^2 \,
 M(x) \, dx + t \, \int_{\R^n} \eta^2(x) \left\vert \nabla u_t(x)
 \right \vert^2 \, M(x) \, dx \\
 = -2 \, t \, \int_{\R^n} \eta(x) \, u_t(x) \, \nabla \eta(x)
 \cdot \nabla u_t(x) \, M(x) \, dx \\
 \leq
 \displaystyle t \, \int_{\R^n} \left\vert u_t(x)\right\vert^2
 \left \vert \nabla \eta(x)\right\vert^2 \, M(x)\, dx
 + t \, \int_{\R^n} \eta^2(x) \left\vert \nabla u_t(x)\right\vert^2
 \, M(x) \, dx,
\end{multline*}
hence
\begin{equation} \label{hence} 
 \int_{\R^n} \eta^2(x)\left\vert
   u_t(x)\right\vert^2 \, M(x) \, dx\leq t \, \int_{\R^n} \left\vert
   u_t(x)\right\vert^2 \left\vert \nabla \eta(x)\right\vert^2 \, M(x)
 \, dx.
\end{equation}
Let $\zeta$ be such that $\zeta=0$ on $E$ and $\zeta$ nonnegative so that $\eta := e^{\alpha \, \zeta}-1 \geq 0$ and $\eta$ vanishes on $E$ for some $\alpha >0$ to be chosen. Choosing this particular $\eta$ in \eqref{hence} with
$\alpha>0$ gives
\[
\int_{\R^n} \left\vert e^{\alpha \, \zeta(x)}-1\right\vert^2\left\vert
 u_t(x)\right\vert^2 \, M(x) \, dx \leq \alpha^2 \, t \, \int_{\R^n}
\left\vert u_t(x)\right\vert^2 \left\vert \nabla \zeta(x)\right\vert^2
\, e^{2 \, \alpha \, \zeta(x)} \, M(x) \, dx.
\]
Taking $\alpha= 1/(2 \, \sqrt{t} \, \left\Vert
 \nabla\zeta\right\Vert_{\infty})$, one obtains
\[
\int_{\R^n} \left\vert e^{\alpha \, \zeta(x)}-1\right\vert^2\left\vert
 u_t(x)\right\vert^2 \, M(x) \, dx\leq \frac 14 \, \int_{\R^n} \left\vert
 u_t(x)\right\vert^2 e^{2\, \alpha \, \zeta(x)} \, M(x) \, dx. 
\]

Using the fact that the norm of $(I+tL)^{-1}$ is bounded by $1$
uniformly in $t >0$, this gives
\[
\begin{array}{lll}
  \displaystyle \left\Vert e^{\alpha\zeta} \, u_t\right\Vert_{L^2(\R^n,M)} 
  & \leq & \displaystyle \left\Vert \left(e^{\alpha\zeta}-1\right) \,  
    u_t\right\Vert_{L^2(\R^n,M)} + \left\Vert u_t\right\Vert_{L^2(\R^n,M)} \\
  & \leq & \displaystyle \frac 12 \left\Vert e^{\alpha\zeta}
    \, u_t\right\Vert_{L^2(\R^n,M)} + \left\Vert f\right\Vert_{L^2(\R^n,M)},
\end{array}
\]
therefore
\[
\begin{array}{lll}
  \displaystyle \int_{\R^n} 
  \left\vert e^{\alpha \, \zeta(x)}\right\vert^2\left\vert
    u_t(x)\right\vert^2 \, M(x) \, dx 
  & \leq & \displaystyle 4 \, \int_{\R^n} \left\vert f(x)\right\vert^2
  \, M(x) \,dx.
\end{array}
\]
We choose now $\zeta$ such that $\zeta =0$ on $E$ as before and
additionnally that $\zeta=1$ on $F$. It can trivially be chosen with
$\left\Vert \nabla\zeta\right\Vert_{\infty} \leq C/d$, which yields the
desired conclusion for the $L^2$ norm of $(I+tL)^{-1}f$ with a factor $4$ in the right-hand side.  Since $t \,
L(\mbox{I}+t \, L)^{-1}f=f-(\mbox{I}+t \, L)^{-1}f$, the desired inequality with a factor $8$
readily follows. \Bk
\hfill\fin
\begin{rem} 
 Arguing similarly, we could also obtain analogous gradient estimates
 for $\left\Vert \sqrt{t} \, \nabla (\mbox{I}+t \, L)^{-1}f\right\Vert_{L^2(F,M)}$. 
\end{rem}
\smallskip

\section{Control of $\left\Vert L^{\alpha/4}f\right\Vert_{L^2(\R^n,M)}$}

This section is devoted to the control of the $L^2$ norm of fractional powers of $L$. This is the 
cornerstone of the proof of Theorem \ref{mainth}.  
In the functional calculus theory of sectorial operators $L$, fractional
powers (for the particular powers we are interested in) are defined as
follows (see for instance \cite[p.24]{henry}:
\begin{equation}\label{def:fractop}
\forall \, \beta \in (0,1), \quad L^{\beta} f =
\frac{1}{\Gamma(1-\beta)} \, \int_0 ^\infty t^{-\beta} \, L \,
e^{-Lt} f \, dt.
\end{equation}
They can also be defined in terms of the resolvent by the Balakrishnan
formulation (see for instance \cite[p.25]{henry}): 
\begin{equation}\label{def:fractop2}
\forall \, \beta \in (0,1), \quad L^{\beta} f=
\frac{\sin(\pi(1-\beta))}{\pi} \, \int_0 ^\infty \lambda^{\beta-1} \,
L \, (L+\lambda \, \mbox{I})^{-1} f \, d\lambda.
\end{equation}

We shall in fact not need any of the representations
\eqref{def:fractop} or \eqref{def:fractop2}; instead we shall rely on the
powerful tool of the so-called ``quadratic estimates'' obtained in the
functional calculus. This is the object of the next lemma. 

\begin{lem} \label{quadratic} Let $\alpha\in (0,2)$. There exists $C_3 =
  C_3(M)>0$ such that, for all $f\in {\mathcal D}(L)$,
\begin{equation} \label{spectral} \left\Vert
   L^{\alpha/4}f\right\Vert_{L^2(\R^n,M)}^2\leq C_3 \,
 \int_0^{+\infty} t^{-1-\alpha/2} \left\Vert t \, L \,
   (\mbox{I} + t \, L)^{-1} f\right\Vert_{L^2(\R^n,M)}^2 \, dt.
\end{equation}
\end{lem}
\noindent {\bf Proof of Lemma \ref{quadratic}.} Let $\mu\in
\left(0,\frac{\pi}2 \right )$ and
$$
\Sigma_{\mu^+}= \left \{ z \in \C^{\ast};\ \left\vert \mbox{arg
   }z\right\vert < \mu\right\} 
$$

$$\Sigma_\mu =\Sigma_{\mu^+} \cup -\Sigma_{\mu^+}. $$

Let $\psi$ be an holomorphic function in $H^\infty(\Sigma_\mu)$ such that
for some $C,\sigma, \tau>0$,
$$
\left\vert \psi(z)\right\vert\leq C \, \inf \left\{ |z|^\sigma \, ; \
  |z|^{-\tau} \right\}
$$
for any $z \in \Sigma_\mu. $ Since $L$ is positive self-adjoint operator on
$L^2(\R^n,M)$ and $L$ is one-to-one on $L^2_0(\R^n,M)$ by (\ref{eq:poincare}), \Bk one has by the spectral theorem
$$
\left\Vert F\right\Vert_{L^2(\R^n,M)}^2\leq C\int_0^{+\infty} \left\Vert
 \psi(tL)F\right\Vert_{L^2(\R^n,M)}^2\frac{dt}t,
$$
whenever $F\in L^2_0(\R^n,M)$. Choosing
$\psi(z)=z^{1-\alpha/4}/(1+z)$ yields
\begin{equation} \label{estimF} \left\Vert F\right\Vert_{L^2(\R^n,M)}^2\leq
 C\int_0^{+\infty} \left\Vert (tL)^{1-\alpha/4} \, (\mbox{I}+t \, L)^{-1}
   \, F\right\Vert_{L^2(\R^n,M)}^2\frac{dt}t
\end{equation}
whenever  $F\in L^2_0(\R^n,M)$. \Bk\par


Let $f\in L^2(\R^n,M)$. Since
$$
\int_{\R^n} Lf(x)M(x)dx=0,
$$
it follows from (\ref{def:fractop}) that the same is true with $L^{\alpha/4}f$. \Bk
Applying now (\ref{estimF}) with $F=L^{\alpha/4}f$ gives the conclusion of
Lemma \ref{quadratic}. \hfill\fin\par
\medskip

Let us draw a simple corollary of Lemma~\ref{quadratic}:
\begin{cor} \label{decoupageent} For any $\alpha \in (0,2)$ and
  $\varepsilon >0$, there is $A = A(M,\varepsilon)$ such that
\begin{equation} \label{decoupe} \left\Vert
   L^{\alpha/4}f\right\Vert_{L^2(\R^n,M)}^2\leq C_3 \,
 \int_0^{A} t^{-1-\alpha/2} \left\Vert t \, L \,
   (\mbox{I} + t \, L)^{-1} f\right\Vert_{L^2(\R^n,M)}^2 \, dt 
 +  \varepsilon \, \| f\|_{L^2(\R^n,M)}^2.
\end{equation}
\end{cor}

\noindent {\bf Proof of Corollary~\ref{decoupageent}.} 
The proof is straightforward since 
$$
\left\Vert t \, L \,
   (\mbox{I} + t \, L)^{-1} f\right\Vert_{L^2(\R^n,M)}^2 \le  C \left\Vert f\right\Vert_{L^2(\R^n,M)}^2 $$
   
  and 
$$   
\int_A ^{+\infty} t^{-1-\alpha/2} \, dt \xrightarrow[A \to +\infty]{}
0.
$$
\hfill\fin\par
\medskip
\Bk

We now come to the desired estimate. 
\begin{lem} \label{controllalpha} Let $\alpha\in (0,2)$ and $\varepsilon$ and $A$ given by Corollary \ref{decoupageent} \Bk. There exist
  $C_4=C_4(M,A)>0$ and $c'=c'(M,A) >0$ such that, for all $f\in {\mathcal
    D}(\R^n)$,
$$
\int_0^{A} t^{-1-\alpha/2} \left\Vert t \, L \, (\mbox{I} + t \, L)^{-1}
  f\right\Vert_{L^2(\R^n,M)}^2 \, dt \Bk \leq 
  $$
  $$C_4 \, \iint_{\R^n \times \R^n}
\frac{\left\vert f(x)-f(y)\right\vert^2}{\left\vert
    x-y\right\vert^{n+\alpha}} \, M(x)\,  e^{-c'\, |x-y|}  \, dx\, dy.
$$
\end{lem} 



\noindent {\bf Proof of Lemma \ref{controllalpha}.} Throughout this proof,
for all $x\in \R^n$ and all $s>0$, denote by $Q(x,s)$ the closed cube
centered at $x$ with side length $s$. For fixed $t \in (0,A)$,
following Lemma~\ref{quadratic}, we shall look for an upper bound for
$\left\Vert t \, L \, (\mbox{I} + t \, L)^{-1} f\right\Vert_{L^2(\R^n,M)}
^2$ involving first order differences for $f$. Pick up a countable family
of points $x_j^t\in \R^n$, $j\in \N$, such that the cubes
$Q\left(x_j^{t},\sqrt{t}\right)$ have pairwise disjoint interiors and
\begin{equation} \label{union}
\R^n=\bigcup_{j \in \N} Q\left(x_j^{t},\sqrt{t}\right).
\end{equation}

By Lemma \ref{cardinal} in Appendix B, there exists a constant $\tilde C>0$
such that for all $\theta>1$ and all $x\in \R^n$, there are at most $\tilde
C\, \theta^n$ indexes $j$ such that $\left\vert x-x_j^t\right\vert\leq
\theta\sqrt{t}$.

\medskip

For fixed $j$, one has
$$
t \, L \, (\mbox{I} + t\, L)^{-1} f= t \, L \, (\mbox{I} + t\, L)^{-1} \,
g^{j,t}
$$
where, for all $x\in \R^n$, 
$$
g^{j,t}(x):=f(x)-m^{j,t}
$$
and $m^{j,t}$ is defined by
$$
m^{j,t}:=\frac 1{\left\vert
   Q\left(x_j^{t},2\sqrt{t}\right)\right\vert}\int_{Q\left(x_j^{t},2\sqrt{t}\right)}
f(y) \, dy
$$
Note that, here, the mean value of $f$ is computed with respect to the
Lebesgue measure on $\R^n$. Since (\ref{union}) holds and the cubes
$Q\left(x_j^{t},\sqrt{t}\right)$ have pairwise disjoint interiors, one
clearly has
$$
\begin{array}{lll}
 \displaystyle \left\Vert t \, L \, (\mbox{I} + t\, L)^{-1} f 
 \right\Vert_{L^2(\R^n,M)}^2 
 & = & \displaystyle \sum_{j \in \N} 
 \left\Vert t \, L \, (\mbox{I} + t\, L)^{-1} f
 \right\Vert_{L^2\left(Q\left(x_j^t,\sqrt{t}\right),M\right)}^2\\
 & = & \displaystyle  \sum_{j \in \N} \left\Vert t\, L \, (\mbox{I} + t\,
   L)^{-1}  g^{j,t}
 \right\Vert_{L^2\left(Q\left(x_j^t,\sqrt{t}\right),M\right)}^2,
\end{array}
$$
and we are left with the task of estimating
$$
\left\Vert t \, L \, (\mbox{I} + t\, L)^{-1}
 g^{j,t}\right\Vert_{L^2\left(Q\left(x_j^{t},\sqrt{t}\right),M\right)} ^2.
$$

To that purpose, set
$$
C_0^{j,t}=Q\left(x_j^t, 2\sqrt{t}\right) \ \mbox{ and
} \ C_k^{j,t}=Q\left(x_j^t,2^{k+1}\sqrt{t}\right)\setminus
Q\left(x_j^t,2^{k}\sqrt{t}\right), \ \forall \, k \ge 1,
$$
and $g^{j,t}_k:=g^{j,t} \, {\bf 1}_{C_k^{j,t}}$, $k \ge 0$, where, for any
subset $A\subset \R^n$, ${\bf 1}_A$ is the usual characteristic function of
$A$. Since $g^{j,t}=\sum_{k\geq 0} g^{j,t}_k$ one has
\begin{equation} 
 \left\Vert t \, L \, (\mbox{I} + t\, L)^{-1}  g^{j,t} 
  \right\Vert_{L^2\left(Q\left(x_j^{t},\sqrt{t}\right),M\right)}
  \leq  \sum_{k\geq 0}  
  \left\Vert t \, L \, (\mbox{I} + t\, L)^{-1}  g_k^{j,t}
   \right\Vert_{L^2\left(Q\left(x_j^{t},\sqrt{t}\right),M\right)} 
\end{equation}
and, using Lemma \ref{off}, one obtains (for some constants $C,c>0$)
\begin{eqnarray} \label{expdecay}
 \left\Vert t \, L \, (\mbox{I} + t\, L)^{-1}  g^{j,t} 
 \right\Vert_{L^2\left(Q\left(x_j^{t},\sqrt{t}\right),M\right)} 
 \leq  
 \\ \nonumber
 C \, \left( \left\Vert
     g_0^{j,t}\right\Vert_{L^2(C_0^{j,t},M)}
   +\sum_{k\geq 1} e^{-c \, 2^{k}} 
   \left\Vert g_k^{j,t}\right\Vert_{L^2(C_k^{j,t},M)} \right).
\end{eqnarray}
By Cauchy-Schwarz's inequality, we deduce (for another constant $C'>0$)
\begin{eqnarray}\label{expdecaybis} 
 \left\Vert t \, L \, (\mbox{I} + t\, L)^{-1}  g^{j,t}
 \right\Vert_{L^2\left(Q\left(x_j^{t},\sqrt{t}\right),M\right)} ^2 \leq \\ \nonumber
 C' \,
 \left( \left\Vert g_0^{j,t}\right\Vert_{L^2(C_0^{j,t},M)}^2 +\sum_{k\geq
     1} e^{-c \, 2^{k}} \left\Vert g_k^{j,t}\right\Vert_{L^2(C_k^{j,t},M)}^2
 \right).
\end{eqnarray}

As a consequence,  we have 
\begin{equation} \label{expdecayter}
\begin{array}{lll}
  \displaystyle \int_0^{A} t^{-1-\alpha/2} 
  \left\Vert t \, L \, (\mbox{I} + t \, L)^{-1}
    f\right\Vert_{L^2(\R^n,M)}^2 \, dt \leq \\ 
  \displaystyle C' \, \int_0^{A}
  t^{-1-\alpha/2} \sum_{j\ge 0} 
  \left\Vert g_0^{j,t}\right\Vert_{L^2(C_0^{j,t},M)}^2 dt+ \\
  \displaystyle C' \, \int_0^{A}  t^{-1-\alpha/2} 
  \sum_{k\geq 1} e^{-c \, 2^{k}} \sum_{j \geq 0}
  \left\Vert g_k^{j,t}\right\Vert_{L^2(C_k^{j,t},M)}^2 dt.
\end{array}
\end{equation}

We claim that:
\begin{lem} \label{estimg}
There exists $\bar C>0$ such that, for all $t>0$ and all $j \in \N$:
\begin{itemize}
\item[{\bf A.}] For the first term:
$$\displaystyle \left\Vert g_0^{j,t}\right\Vert_{L^2(C_0^{j,t},M)}^2\leq  
\frac{\bar C}{t^{n/2}} \int_{Q\left(x_j^t,2\sqrt{t}\right)}
\int_{Q\left(x_j^t,2\sqrt{t}\right)} \left\vert f(x)-f(y)\right\vert^2 \,
M(x) \, dx \, dy.$$
\item[{\bf B.}]
For all $k\geq 1$,
\[ \left\Vert g^{j,t}_k\right\Vert_{L^2(C_k^{j,t},M)}^2 
\leq \]
\[
 \frac{\bar C}{(2^k\sqrt{t})^n} \int_{x\in Q(x^t_j,2^{k+1}\sqrt{t})} 
\int_{y\in Q(x^t_j,2^{k+1}\sqrt{t})} \left\vert f(x)-f(y)\right\vert^2 \, M(x)
\, dx \, dy.\]
\end{itemize}
\end{lem}

We postpone the proof to the end of the section and finish the proof of
Lemma \ref{controllalpha}. Using Assertion {\bf A} in Lemma \ref{estimg},
summing up on $j \ge 0$ and integrating over $(0,A)$, we get 
\begin{multline*}
 \displaystyle \int_0^{A} t^{-1-\alpha/2} \sum_{j \ge 0} \left \Vert
   g_0^{j,t}\right\Vert_{L^2\left(C_0^{j,t},M\right)}^2 \, dt = \sum_{j
   \ge 0} \int_0^{A} t^{-1-\alpha/2} \left\Vert g_0^{j,t}\right
 \Vert_{L^2\left(C_0^{j,t},M\right)}^2 \, dt \\
 \displaystyle \le \bar C \, \sum_{j \ge 0} \int_0^{A}
 t^{-1-\frac{\alpha}2-\frac n2} \left(\int_{Q\left(x_j^t,2\sqrt{t}\right)}
   \int_{Q\left(x_j^t,2\sqrt{t}\right)} \left\vert
     f(x)-f(y)\right\vert^2 \, M(x) \, dx \, dy\right) \, dt \\
 \displaystyle \le \bar C\, \sum_{j \ge 0}\iint_{(x,y)\in \R^n\times \R^n}
 \left\vert f(x)-f(y)\right\vert^2M(x) \times \\
 \left(\int_{t\geq
     \max\left\{\frac{\left\vert x-x_j^t\right\vert^2}n\,;\
       \frac{\left\vert y-x_j^t\right\vert^2}n\right\}} ^A
   \, t^{-1-\frac{\alpha}2-\frac n2}dt\right) \, dx \, dy.
\end{multline*}

The Fubini theorem now shows
$$
\sum_{j \ge 0} \int_{t\geq \max\left\{\frac{\left\vert
       x-x_j^t\right\vert^2}n\, ; \ \frac{\left\vert
       y-x_j^t\right\vert^2}n \right\}} ^A \, t^{-1-\frac{\alpha}2-\frac n2} \, dt
= $$

$$
\int_0^{A} t^{-1-\frac{\alpha}2-\frac n2} \, \sum_{j \ge 0} {\bf
 1}_{\left(\max\left\{ \frac{\left\vert x-x_j^t\right\vert^2}n\, ; \
     \frac{\left\vert y-x_j^t\right\vert^2}n\right\},+\infty\right)}
(t) \, dt.
$$
Observe that, by Lemma \ref{cardinal}, there is a constant $N \in \N$ such
that, for all $t>0$, there are at most $N$ indexes $j$ such that
$\left\vert x-x_j^t\right\vert^2< n\, t$ and $\left\vert
  y-x_j^t\right\vert^2<n \, t$, and for these indexes $j$, one has
$\left\vert x-y\right\vert<2 \, \sqrt{nt}$. It therefore follows that
$$
\sum_{j \ge 0} {\bf 1}_{\left(\max\left\{ \frac{\left\vert
          x-x_j^t\right\vert^2}n \, ; \ \frac{\left\vert
          y-x_j^t\right\vert^2}n\right\},+\infty\right)}(t)\leq N \, {\bf
  1}_{\left(\left\vert x-y\right\vert^2/4n,+\infty\right)}(t),
$$
so that
\begin{multline} \label{intg0}
 \displaystyle \int_0^{A} t^{-1-\alpha/2} 
 \sum_{j} \left\Vert
   g_0^{j,t}\right\Vert_{L^2\left(C_0^{j,t},M\right)}^2 \, dt \\
 \leq \bar C \, N \, \iint_{\R^n\times \R^n} 
  \left\vert f(x)-f(y)\right\vert^2M(x) 
  \left( \int_{\left\vert x-y\right\vert^2/4n}^{A} \,
   t^{-1-\frac{\alpha}2-\frac{n}2} \, dt\right) \, dx \, dy \Bk\\
 \displaystyle  \leq  \bar C \, N \, \iint_{|x-y|\le 2\sqrt{nA}\Bk} 
  \frac{\left\vert f(x)-f(y)\right\vert^2}
   {\left\vert x-y\right\vert^{n+\alpha}} \, M(x) \, dx \, dy.
\end{multline}
Using now Assertion {\bf B} in Lemma \ref{estimg}, we obtain, for all $j
\ge 0$ and all $k\geq 1$,

$$
\begin{array}{l}
  \displaystyle \int_0^{A} t^{-1-\alpha/2} 
  \, \sum_{j \ge 0} \left\Vert g^{j,t}_k\right\Vert_2^2dt  \\
  \displaystyle \leq 
  \bar C \, 2^{-kn} \, \sum_{j \ge 0} 
  \int_0^{A} t^{-1-\frac{\alpha}2-\frac{n}2} \, 
  \left(\iint_{Q(x^t_j,2^{k+1}\sqrt{t}) \times  
      Q(x^t_j,2^{k+1}\sqrt{t})} \left\vert f(x)-f(y)\right\vert^2 \, M(x) \, dx
    \, dy\right) \, dt\\
  \displaystyle \leq \bar C \, 2^{-kn} \, 
  \sum_{j \ge 0} \iint_{x,y\in \R^n} 
  \left\vert f(x)-f(y)\right\vert^2 \, M(x) \times \\ \, 
  \left(\int_0^{A} t^{-1-\frac{\alpha}2-\frac{n}2} \, 
    {\bf 1}_{\left(\max\left\{\frac{\left\vert x-x^t_j\right\vert^2}{4^{k}n}, 
   \frac{\left\vert y-x^t_j\right\vert^2}{4^{k}n}\right\},+\infty\right)} 
    (t) \, dt\right) \, dx \, dy. 
\end{array}
$$
But, given $t>0$, $x,y\in \R^n$, by Lemma \ref{cardinal} again, there exist
at most $\tilde C \, 2^{kn}$ indexes $j$ such that
$$
\left\vert x-x_j^t\right\vert\leq 2^{k}\sqrt{nt} \ \mbox{ and } \ 
\left\vert y-x_j^t\right\vert\leq 2^{k}\sqrt{nt}, 
$$
and for these indexes $j$, $\left\vert x-y\right\vert\leq
2^{k+1}\sqrt{nt}$. As a consequence,
\begin{equation} \label{intgk}
\begin{array}{lll}
  \displaystyle 
  \int_0^{A} t^{-1-\frac{\alpha}2-\frac{n}2} \, 
  \sum_{j \ge 0} {\bf 1}_{\left(\max\left\{\frac{\left\vert
            x-x^t_j\right\vert^2}{4^{k}n}, 
  \frac{\left\vert x-x^t_j\right\vert^2}{4^{k}n}\right\},+\infty\right)}(t)
  \, dt  \leq \\
  \displaystyle \tilde C \, 2^{kn} \, 
  \int_{t\geq \frac{\left\vert x-y\right\vert^2}{4^{k+1}n}} ^{A}  \,
  t^{-1-\frac{\alpha}2-\frac{n}2} \, dt \leq \\
  \displaystyle \tilde C' \, 2^{k(\alpha+n)} \, 
  \left\vert x-y\right\vert^{-n-\alpha} \, {\bf 1}_{|x-y|\le 2^{k+1} \,
    \sqrt{n \, A}},
\end{array}
\end{equation}
for some other constant $\tilde C' >0$, and therefore
$$
 \displaystyle \int_0^{A} t^{-1-\alpha/2} 
 \sum_{j} \left\Vert
   g_k ^{j,t}\right\Vert_{L^2\left(C_0^{j,t},M\right)}^2 \, dt \leq 
 $$
 
 $$
   \bar C \, \tilde C' \, 2^{k(\alpha+n)} \, \iint_{|x-y|\le 2^{k+1} \,
    \sqrt{n \, A}} 
  \frac{\left\vert f(x)-f(y)\right\vert^2}
   {\left\vert x-y\right\vert^{n+\alpha}} \, M(x) \, dx \, dy.
$$

We can now conclude the proof of Lemma \ref{controllalpha}, using
Lemma \ref{quadratic}, (\ref{expdecay}), (\ref{intg0}) and (\ref{intgk}). 
We have proved, by reconsidering \eqref{expdecayter}: 
\begin{multline}
\displaystyle 
\int_0^{A} t^{-1-\alpha/2} \left\Vert t \, L \, (\mbox{I} + t \, L)^{-1}
  f\right\Vert_{L^2(\R^n,M)}^2 \, dt
\leq 
\\
\displaystyle C' \, \bar C \, N \, \iint_{|x-y|\le  2 \,
  \sqrt{n \, A}\Bk} 
\frac{\left\vert f(x)-f(y)\right\vert^2}
{\left\vert x-y\right\vert^{n+\alpha}} \, M(x) \, dx \, dy \\
+ \displaystyle \sum_{k \ge 1} 
C' \, \bar C \, \tilde C' \, 2^{k\alpha} \, e^{-c \, 2^k} \,  
\iint_{|x-y|\le 2^{k+1} \, \sqrt{n \, A}} 
\frac{\left\vert f(x)-f(y)\right\vert^2}
{\left\vert x-y\right\vert^{n+\alpha}} \, M(x) \, dx \, dy
\end{multline}
and we deduce that 
$$
\displaystyle 
\int_0^{A} t^{-1-\alpha/2} \left\Vert t \, L \, (\mbox{I} + t \, L)^{-1}
  f\right\Vert_{L^2(\R^n,M)}^2 \, dt
\le $$
$$
C_4 \, \iint_{\R^n\times \R^n} \frac{\left\vert f(x)-f(y)\right\vert^2}
{\left\vert x-y\right\vert^{n+\alpha}} \, M(x) \, e^{-c' \, |x-y|}
\, dx \, dy
$$
for some constants $C_4$ and $c'>0$ as claimed in the statement.
\Bk \hfill\fin\par

\medskip

\noindent
{\bf Proof of Lemma \ref{estimg}.} Observe first that, for all $x\in \R^n$,
$$
\begin{array}{lll}
 \displaystyle g_0^{j,t}(x) & = & \displaystyle f(x) 
  -\frac 1{\left\vert Q\left(x_j^t,2\sqrt{t}\right)\right\vert} 
   \int_{Q\left(x_j^t,2\sqrt{t}\right)} f(y) \, dy\\
 & = & \displaystyle \frac 1{\left\vert 
  Q\left(x_j^t,2\sqrt{t}\right)\right\vert} 
   \int_{Q\left(x_j^t,2\sqrt{t}\right)} (f(x)-f(y)) \, dy.
\end{array}
$$
By Cauchy-Schwarz inequality, it follows that
$$
\left\vert g_0^{j,t}(x)\right\vert^2\leq \frac C{t^{n/2}}   \int_{Q\left(x_j^t,2\sqrt{t}\right)} \left\vert f(x)-f(y)\right\vert^2dy.
$$
Therefore,
$$
\left\Vert g_0^{j,t}\right\Vert_{L^2(C_0^{j,t},M)}^2\leq \frac
C{t^{n/2}} \int_{Q\left(x_j^t,2\sqrt{t}\right)}
\int_{Q\left(x_j^t,2\sqrt{t}\right)} \left\vert
 f(x)-f(y)\right\vert^2 \, M(x) \, dx \, dy,
$$
which shows Assertion {\bf A}. We argue similarly for Assertion {\bf B} and
obtain
$$
\displaystyle \left\Vert g^{j,t}_k\right\Vert_{L^2 (C_k ^{j,t},M)}^2 
\leq \displaystyle \frac C{2^{k/n} t^{n/2}} \int_{x\in Q(x^t_j,2^{k+1}\sqrt{t})}
\int_{y\in Q(x^t_j,2^{k+1}\sqrt{t})} \left\vert
 f(x)-f(y)\right\vert^2 \, M(x) \, dx \, dy,
$$
which ends the proof of Lemma \ref{estimg}. \hfill\fin\par

\bigskip

We end up this section with a few comments on Lemma \ref{estimg}. It is a
well-known fact (\cite{stri}) that, when $0<\alpha<2$, for all $p\in
(1,+\infty)$,
\begin{equation} \label{compars} \left\Vert
    (-\Delta)^{\alpha/4}f\right\Vert_{L^p(\R^n)} \leq C_{\alpha,p}
  \left\Vert S_{\alpha}f\right\Vert_{L^p(\R^n)}
\end{equation}
where
$$
S_{\alpha}f(x)=\left(\int_0^{+\infty} \left(\int_B \left\vert
      f(x+ry)-f(x)\right\vert dy\right)^2
  \frac{dr}{r^{1+\alpha}}\right)^{\frac 12},
$$
and also (\cite{stein}) 
\begin{equation} \label{compard} \left\Vert
    (-\Delta)^{\alpha/4}f\right\Vert_{L^p(\R^n)} \leq C_{\alpha,p}
  \left\Vert D_{\alpha}f\right\Vert_{L^p(\R^n)}
\end{equation}
where
$$
D_{\alpha}f(x)=\left(\int_{\R^n} \frac{\left\vert
      f(x+y)-f(x)\right\vert^2}{\left\vert
      y\right\vert^{n+\alpha}}dy\right)^{\frac 12}.
$$
In \cite{crt}, these inequalities were extended to the setting of a
unimodular Lie group endowed with a sub-laplacian $\Delta$, relying on
semigroups techniques and Littlewood-Paley-Stein functionals. In
particular, in \cite{crt}, we use {\it pointwise} estimates of the kernel
of the semigroup generated by $\Delta$. The conclusion of Lemma
\ref{estimg} means that the norm of $L^{\alpha/4}f$ in $L^2(\R^n,M)$ is
bounded from above by the $L^2(\R^n,M)$ norm of an appropriate version of
$D_{\alpha}$. Note that this does not require pointwise estimates for the
kernel of the semigroup generated by $L$, and that the $L^2$ off-diagonal
estimates given by Lemma \ref{off}, which hold for a general measure $M$,
are enough for our argument to hold. However, we do not know if an $L^p$
version of Lemma \ref{estimg} still holds. Note also that we do not compare
the $L^2(\R^n,M)$ norm of $L^{\alpha/4}f$ with the $L^2(\R^n,M)$ norm of a
version of $S_{\alpha}f$. Finally, the converse inequalities to
(\ref{compars}) and (\ref{compard}) hold in $\R^n$ and also on a unimodular
Lie group (\cite{crt}), and we did not consider the corresponding
inequalities in the present paper. \Bk

\section{Control of the moment of $f$ by $\left\Vert
    L^{\alpha/4}f\right\Vert_{L^2(\R^n,M)}$ and proof of Theorem
  \ref{mainth}} \label{sec:fin}

Observe first that, by the definition of $L$, we have
$$
\int_{\R^n} \left\vert \nabla f(x)\right\vert^2 \, M(x) \, dx=\int_{\R^n}
Lf(x) \, f(x) \, M(x) \, dx.
$$
for all $f\in {\mathcal D}(L)$. The inequality (\ref{pim}) can therefore be
rewritten, in terms of operators, as
\begin{equation} \label{ineqop}
L\geq \lambda' \, \mu, 
\end{equation}
where $\mu$ is the multiplication operator by $x\mapsto 1+\left\vert \nabla
  \ln M(x)\right\vert^2$. Since $\mu$ is a nonnegative operator on
$L^2(\R^n,M)$, using a functional calculus argument (see \cite{davies}, p.
110), one deduces from (\ref{ineqop}) that, for any $\alpha\in (0,2)$,
\[
L^{\alpha/2}\geq (\lambda')^{\alpha/2} \, \mu^{\alpha/2},
\]
which implies, thanks to the fact $L^{\alpha/2} = (L^{\alpha/4})^2$ and the
symmetry of $L^{\alpha/4}$ on $L^2(\R^n,M)$, that
$$
(\lambda')^{\alpha/2} \, \int_{\R^n} \left\vert f(x)\right\vert^2
(1+\left\vert \nabla\ln M(x)\right\vert^2)^{\alpha/2} \, M(x) \, dx\leq $$
$$\int_{\R^n} \left\vert L^{\alpha/4}f(x)\right\vert^2 \, M(x) \, dx =
\left\| L^{\alpha/4} f \right\|^2 _{L^2(\R^n,M)}.
$$
The conclusion of Theorem \ref{mainth} readily follows by using the
previous inequality in conjunction with Lemma \ref{decoupageent} and
\ref{controllalpha}, and picking $\varepsilon$ small enough.

\section{Appendix A: Improved Poincar\'e inequality}

In this section, we prove Proposition \ref{poincimprovedM}, namely

\begin{pro}  
 Assume that $M=e^{-V}$ satisfies (\ref{alpha}). 
Then there exists $\lambda'(M)>0$ such that, for
 all function $f\in L^2_0(\R^n,M) \cap H^1(\R^n,M)$:
 \begin{equation} 
   \int_{\R^n} \left\vert \nabla f(x)\right\vert^2 \, M(x) \, dx 
   \ge \lambda'(M) \, \int_{\R^n} \left\vert f(x)\right\vert
   ^2\left(1+\left\vert \nabla \ln M(x)\right\vert^2\right) \, M(x) \, dx.
\end{equation}
\end{pro}

Note that
of course in general the constants $\lambda(M)$ and $\lambda'(M)$ in
\eqref{eq:poincare} and \eqref{pim} are different.\par

\noindent
{\bf Proof of Proposition \ref{poincimprovedM}.} Let $f$ be as in the
statement of Proposition \ref{poincimprovedM} and let $g:=fM^{\frac 12}$.
Since
$$
\nabla f=M^{-\frac 12} \, \nabla g-\frac 12 \, g \, M^{-\frac 32} \, \nabla M,
$$
assumption (\ref{alpha}) yields two positive constants $\beta,\gamma$ such that
\begin{multline} \label{estim1}
 \displaystyle \int_{\R^n} \left\vert \nabla f(x)\right\vert^2(x) \,
 M(x) \, dx  = \\
\displaystyle \int_{\R^n} \left(\left\vert \nabla g(x)\right\vert^2 
+\frac 14g^2(x)\left\vert \nabla \ln M(x)\right\vert^2- 
g(x)\nabla g(x)\cdot \nabla \ln M(x)\right) \, dx\\
 = \displaystyle \int_{\R^n} \left(\left\vert \nabla
     g(x)\right\vert^2 
+\frac 14g^2(x)\left\vert \nabla \ln M(x)\right\vert^2-\frac 12
\nabla g^2(x)\cdot \nabla \ln M(x)\right) \, dx\\
 \geq \displaystyle \int_{\R^n} g^2(x) \left(\frac 14\left\vert
     \nabla \ln M(x)\right\vert^2+\frac 12 \Delta\ln M(x)\right) \, dx\\
 \geq  \displaystyle \int_{\R^n} f^2(x)\left(\beta \left\vert
     \nabla \ln M(x)\right\vert^2-\gamma\right) \, M(x) \, dx.
\end{multline}
The conjunction of \eqref{eq:poincare} (which holds because (\ref{assumptionV}) is satisfied),  and (\ref{estim1}) yields the
desired conclusion. \hfill\fin\par

\section{Appendix B: Technical lemma}

We prove the following lemma. 
\begin{lem} \label{cardinal} There exists a constant $\tilde C>0$ with the
 following property: for all $\theta>1$ and all $x\in \R^n$, there are at
 most $\tilde C\, \theta^n$ indexes $j$ such that $\left\vert
   x-x_j^t\right\vert\leq \theta\sqrt{t}$.
\end{lem} 

\noindent {\bf Proof of Lemma~\ref{cardinal}.} The argument is very simple
(see \cite{kanai}) and we give it for the sake of completeness. Let $x\in \R^n$
and $I(x):= \left\{j \in \N \, ;\ \left\vert x-x_j^t\right\vert\leq
 \theta\sqrt{t}\right\}$.  Since, for all $j\in I(x)$,
$$Q\left(x_j^t,\sqrt{t}\right)\subset B\left(x,\left(\theta+\frac
   12\right)\sqrt{nt}\right),$$ one has
$$
C\left(\left(\theta+\frac 12\right)\sqrt{nt}\right)^n\geq \sum_{j\in I(x)}
\left\vert Q\left(x_j^t,\sqrt{t}\right)\right\vert= |I(x)| \, 
\sqrt{t}^n,
$$
we get the desired conclusion. \hfill\fin\par

\bigskip

\noindent {\bf Acknowledgement}: The first author would like to thank
the Award No. KUK-I1-007-43, funded by the King Abdullah University of
Science and Technology (KAUST) for the funding provided in Cambridge
University.  \bigskip

 \bibliographystyle{alpha} 
  \bibliography{Biblio-puiss-frac}

\newcommand{\etalchar}[1]{$^{#1}$}
\begin{thebibliography}{CRTN01}

\bibitem[AF03]{adams}
R.A. Adams and J.J.F. Fournier.
\newblock {\em Sobolev spaces}, volume 140 of {\em Pure and Applied Mathematics
  (Amsterdam)}.
\newblock Elsevier/Academic Press, Amsterdam, second edition, 2003.

\bibitem[AHL{\etalchar{+}}02]{kato}
P.~Auscher, S.~Hofmann, M.~Lacey, A.~McIntosh, and P.~Tchamitchian.
\newblock {The solution of the Kato square root problem for second order
  elliptic operators on ${\mathbb R}^n$}.
\newblock {\em Ann. of Math. (2)}, 156(2):633--654, 2002.

\bibitem[AMR08]{amr}
P.~Auscher, A.~McIntosh, and E.~Russ.
\newblock {Hardy spaces of differential forms on Riemannian manifolds}.
\newblock {\em J. Geom. Anal.}, 18(1):192--248, 2008.

\bibitem[Aus07]{auscherams}
P.~Auscher.
\newblock {\em {On necessary and sufficient conditions for $L^p$ estimates of
  Riesz transforms associated to elliptic operators on $\R^n$ and related
  estimates}}, volume 186 of {\em Mem. Amer. Math. Soc}.
\newblock Amer. Math. Soc., 2007.

\bibitem[BBCG08]{BBCG}
Dominique Bakry, Franck Barthe, Patrick Cattiaux, and Arnaud Guillin.
\newblock A simple proof of the {P}oincar\'e inequality for a large class of
  probability measures including the log-concave case.
\newblock {\em Electron. Commun. Probab.}, 13:60--66, 2008.

\bibitem[B{\'E}86]{be}
D.~Bakry and M.~{\'E}mery.
\newblock Propaganda for {$\Gamma\sb 2$}.
\newblock In {\em From local times to global geometry, control and physics
  ({C}oventry, 1984/85)}, volume 150 of {\em Pitman Res. Notes Math. Ser.},
  pages 39--46. Longman Sci. Tech., Harlow, 1986.

\bibitem[Cha04]{chafai}
D.~Chafa{\i}.
\newblock {Entropies, convexity, and functional inequalities: on
  $\Phi$-entropies and $\Phi$-Sobolev inequalities}.
\newblock {\em J. Math. Kyoto Univ.}, 44(2):325--363, 2004.

\bibitem[CRTN01]{crt}
T.~Coulhon, E.~Russ, and V.~Tardivel-Nachef.
\newblock {Sobolev algebras on Lie groups and Riemannian manifolds}.
\newblock {\em Amer. J. Math.}, 123:283--342, 2001.

\bibitem[Dav80]{davies}
E.B. Davies.
\newblock {\em One-parameter semigroups}, volume~15 of {\em London Mathematical
  Society Monographs}.
\newblock Academic Press, Inc., London-New York, second edition, 1980.

\bibitem[DS90]{deuschstroock}
Jean-Dominique Deuschel and Daniel~W. Stroock.
\newblock Hypercontractivity and spectral gap of symmetric diffusions with
  applications to the stochastic {I}sing models.
\newblock {\em J. Funct. Anal.}, 92(1):30--48, 1990.

\bibitem[Gaf59]{gaff}
M.P. Gaffney.
\newblock The conservation property of the heat equation on {R}iemannian
  manifolds.
\newblock {\em Comm. Pure Appl. Math.}, 12:1--11, 1959.

\bibitem[GI08]{GI}
I.~Gentil and C.~Imbert.
\newblock {The L{\'e}vy--Fokker--Planck equation: $\Phi$-entropies and
  convergence to equilibrium}.
\newblock {\em Asymptot. Anal.}, 59(3):125--138, 2008.

\bibitem[Hen81]{henry}
D.~Henry.
\newblock {\em Geometric theory of semilinear parabolic equations}, volume 840
  of {\em Lecture Notes in Mathematics}.
\newblock Springer-Verlag, Berlin-New York, second edition, 1981.

\bibitem[Kan85]{kanai}
M.~Kanai.
\newblock Rough isometries and combinatorial approximations of geometries of
  noncompact riemannian manifolds.
\newblock {\em J. Math. Soc. Japan}, 37:391--413, 1985.

\bibitem[Lan72]{landkof}
N.~S. Landkof.
\newblock {\em Foundations of modern potential theory}.
\newblock Springer-Verlag, New York, 1972.
\newblock Translated from the Russian by A. P. Doohovskoy, Die Grundlehren der
  mathematischen Wissenschaften, Band 180.

\bibitem[Led01]{ledouxlivre}
M.~Ledoux.
\newblock {\em {The concentration of measure phenomenon}}.
\newblock Amer. Math. Soc., 2001.

\bibitem[Ste61]{stein}
E.M. Stein.
\newblock {The characterization of functions arising as potentials I}.
\newblock {\em Bull. Amer. Math. Soc.}, 67:102--104, 1961.

\bibitem[Str67]{stri}
R.S. Strichartz.
\newblock Multipliers on fractional {S}obolev spaces.
\newblock {\em J. Math. Mech.}, 16:1031--1060, 1967.

\bibitem[Vil09]{villaniams}
C.~Villani.
\newblock {HypocoercivityI}.
\newblock {\em To appear in Memoirs of the AMS}, 2009.

\bibitem[Wu00]{wu}
L.~Wu.
\newblock {A new modified logarithmic Sobolev inequality for Poisson point
  processes and several applications}.
\newblock {\em Probab. Theory Related Fields}, 118(3):427--438, 2000.

\end{thebibliography}
\bigskip

\noindent      
{\em Cl\'ement Mouhot}-- University of Cambridge, DAMTP\\
Wilberforce road, Cambridge CB3 0WA, England\\
{\it On leave from:} CNRS \& \'Ecole Normale
Sup\'erieure, \\
DMA, 45, rue d'Ulm - F 75230 Paris cedex 05, France
   
\medskip 
   
\noindent
{\em Emmanuel Russ}--
Universit\'e Aix-Marseille III, LATP,\\
Facult\'e des Sciences et Techniques, Case cour A\\
Avenue Escadrille Normandie-Niemen, F-13397 Marseille,
Cedex 20, France et  \\
CNRS, LATP, CMI, 39 rueF. Joliot-Curie, F-13453 Marseille Cedex 13,
France

\medskip 

\noindent
{\em Yannick Sire}--
Universit\'e Aix-Marseille III, LATP,\\
Facult\'e des Sciences et Techniques, Case cour A\\
Avenue Escadrille Normandie-Niemen, F-13397 Marseille,
Cedex 20, France et  \\
CNRS, LATP, CMI, 39 rueF. Joliot-Curie, F-13453 Marseille Cedex 13,
France

\end{document}